\documentclass[a4paper, 12pt]{article}
\usepackage{epsfig, color, epstopdf}
\usepackage{amssymb,amsmath,amsthm,latexsym,lineno}
\usepackage{graphicx,graphics,amsfonts,indentfirst}
\usepackage{mathrsfs,eufrak, bbm, dsfont,yfonts}
\usepackage{url,tikz,yhmath,appendix,booktabs}
\usepackage{multirow}
\usepackage{kantlipsum}
\allowdisplaybreaks
\usepackage{setspace}
\usepackage{authblk}

\newtheorem{theo}{Theorem}[section]
\newtheorem{lem}[theo]{Lemma}
\newtheorem{pro}[theo]{Proposition}
\newtheorem{prob}{Problem}
\newtheorem{exa}[theo]{Example}
\newtheorem{con}{Conjecture}
\newtheorem{cor}[theo]{Corollary}

\newtheorem{obs}{Observation}

\usetikzlibrary{backgrounds}
\usetikzlibrary{topaths,calc}
\usetikzlibrary{arrows}
\usetikzlibrary{shapes,shapes.geometric,shapes.misc}
\pgfdeclarelayer{edgelayer}
\pgfdeclarelayer{nodelayer}
\pgfsetlayers{background,edgelayer,nodelayer,main}
\tikzstyle{none}=[inner sep=0mm]
\tikzstyle{every loop}=[]

\tikzstyle{dotted}=[dash pattern=on \pgflinewidth off 2pt]
\tikzstyle{dashed}=[dash pattern=on 3pt off 3pt]

\tikzstyle{new style 0}=[fill=black, draw=black, shape=circle]
\tikzstyle{red style 1}=[fill=red, draw=black, shape=circle]
\tikzstyle{blue style 2}=[fill=blue, draw=black, shape=circle]
\tikzstyle{white style 4}=[fill=white, draw=black, shape=circle]
\tikzstyle{bklack style 5}=[fill=black, draw=black, shape=rectangle]
\tikzstyle{red style 3}=[fill=red, draw=black, shape=rectangle]
\tikzstyle{yellow style 7}=[fill=yellow, draw=black, shape=rectangle]
\tikzstyle{new style 8}=[fill={rgb,255: red,0; green,132; blue,0}, draw={rgb,255: red,0; green,131; blue,0}, shape=circle]

\tikzstyle{new edge style 0}=[-]
\tikzstyle{new edge style 1}=[-, draw=red]
\tikzstyle{new edge style 2}=[-, draw=blue]
\tikzstyle{new edge style 3}=[-, draw={rgb,255: red,0; green,156; blue,0}]

\tikzstyle{cblue}=[circle, draw, thin,fill=blue!20, scale=0.5]

\newcommand \equ[2]
{
	\begin{equation}\label{#1}
		#2
	\end{equation}
}

\newcommand \aln[2]
{
	\begin{align}\label{#1}
		#2
	\end{align}
}

\newcommand \eqn[2]
{
	\begin{eqnarray}\label{#1}
		#2
	\end{eqnarray}
}

\newcommand \lemm[2]
{\begin{lem}
		\label{#1} #2
	\end{lem}
}

\newcommand \prop[2]
{\begin{pro}
		\label{#1} #2
	\end{pro}
}

\newcommand \corr[2]
{\begin{cor}
		\label{#1} #2
	\end{cor}
}

\def \setd {{\cal D}}

\newcounter{countcase}

\newcounter{countclaim}

\def \proof {\noindent {{\it Proof}}.\setcounter{countcase}{0} \setcounter{clm}{0}}

\newcommand{\proofend}{{\hfill$\Box$}\setcounter{countclaim}{0}\setcounter{countcase}{0}}

\def \N {{\mathbb N}}

\def \hyh {{\cal H}}
\newcommand \spann[1]{\langle #1\rangle }

\def \nbroken {{\mathscr {N}}\hspace{-0.15 cm}{\mathscr {B}}}

\usepackage{xpatch}
\makeatletter
\xpatchcmd{\@thm}{\thm@headpunct{.}}{\thm@headpunct{}}{}{}
\makeatother

\begin{document}
	\baselineskip 0.6 cm
	
	\title{
		Compare list-color functions
		of uniform hypergraphs with their chromatic polynomials
	}

	\author[1]{\small Fengming Dong\thanks{Corresponding Author. Email: fengming.dong@nie.edu.sg}}
		\author[2]{\small Meiqiao Zhang\thanks{Email: 
			meiqiaozhang95@163.com}}

	\affil[1]{\footnotesize 
		National Institute of Education,
		Nanyang Technological University, 
		Singapore}
	
		\affil[2]{\footnotesize School of Mathematical Sciences, Xiamen  University, China}
	
	\date{}
	
	\maketitle{}

	\begin{abstract}
		
	In [J. Combin. Theory Ser. B 161
		(2023), 109--119], the authors 
		showed that the list-color function
		$P_l(G,k)$ of any simple graph $G$ of size $m$ coincides 
		with its chromatic polynomial $P(G,k)$ for 
		all integers $k\ge m-1$. 
		In this article, we extend this  conclusion to 
		any uniform hypergraph. 
		Furthermore, we show that 
		for any $r$-uniform hypergraph 
		${\cal H}=(V,E)$,
		 where $r\ge 2$,
		$P({\cal H}, L)-P({\cal H},k)\ge 
		(k-|E|+1)k^{|V|-r-1}
		\sum\limits_{e\in E}
		\left (k-\left|\bigcap\limits_{v\in e}L(v)\right|\right ) 
		$ holds for all integers $k$ with $k\ge |E|-1\ge 4$,
		where $L$ is any $k$-assignment of ${\cal H}$ 
		and $P({\cal H}, L)$ 
		is the number of $L$-colorings of 
		${\cal H}$. 
	\end{abstract}


	\smallskip
	


\section{ Introduction
\label{secintro}}

In this article, we will consider hypergraphs in which each edge has at least two vertices and 
no edge is contained in another edge. 
For a hypergraph $\hyh$, 
let $V(\hyh)$ and $E(\hyh)$ denote its
vertex set and edge set respectively. 
By convention, we say $\hyh$ is \textit{$r$-uniform} if $|e|=r$ for each edge $e\in E(\hyh)$.
Let $c(\hyh)$ be the number of connected components of $\hyh$.
For any non-empty subset 
$S\subseteq V(\hyh)$, 
let $\hyh[S]$ denote the subhypergraph of $\hyh$ induced by $S$.
For any edge set $F\subseteq E(\hyh)$, 
let $\hyh\spann{F}$ be the spanning subhypergraph of $\hyh$ with edge set $F$, let $c(F)=c(\hyh\spann{F})$, and
let $V(F)$ be the set of vertices 
$v\in V(\hyh)$ with $v\in e$ for some $e\in F$.
Thus, $V(F)=\emptyset$ if and only if 
$F=\emptyset$.
When $F\ne \emptyset$, 
let $\hyh[F]$ be the subhypergraph of $\hyh$ induced by $F$,
i.e., the subhypergraph
with vertex set $V(F)$ and edge set $F$.
For any real number $w$ and positive integer $s$, write $w_{(s)}$ 
for the expression $w(w-1)\cdots (w-s+1)$.

\subsection{Chromatic polynomials	and list-color functions
}

Let $\N$ denote the set of positive integers and $[k]=\{1,2,\dots,k\}$ for any $k\in\N$.
For any $k\in\N$, a \textit{proper $k$-coloring} of a graph $G$ 
(resp., a hypergraph $\hyh$) is a mapping $\theta: V(G)\rightarrow [k]$
(resp. $\theta:V(\hyh)\rightarrow [k]$)
such that $\theta(u)\ne \theta(v)$ 
for each edge $e=uv\in E(G)$
(resp., $|\{\theta(v):   v\in e\}|\ge 2$ for each $e\in E(\hyh)$).
The \textit{chromatic polynomial} 
$P(G, k)$ of $G$,  introduced by Birkhoff~\cite{birk} in 1912, is the polynomial  counting the number of proper $k$-colorings of $G$, for each $k\in\N$. 
Although the original purpose of applying 
the chromatic polynomial 
to attack the four-color conjecture 
has not been successful yet, 
many of its elegant properties and various connections to other graph-polynomials
were discovered along the way, and it has thus become an interesting research topic. 
See~\cite{dong1,dong0,Jack15,rea1,rea2, Whitney1932} for further information.
Then naturally, the chromatic polynomial $P(\hyh, k)$ of a hypergraph $\hyh$ is defined to count the number of proper $k$-colorings of $\hyh$ for any $k\in\N$. 
In the literature, whether certain well-known properties of the chromatic polynomials of graphs hold for hypergraphs or not has always been a research focus. 
Some recent works in this area 
are referred to~\cite{Dur2022,  
	Tomescu1998, Tomescu2009, Trinks14,
	wang20, ruixue1,ruixue2,ruixue3}.

{\it List-coloring}, introduced independently by 
 Vizing~\cite{vizing} and Erd\H{o}s, Rubin and Taylor~\cite{erdos},
 is a generalization of proper coloring.
For any $k\in\N$, a \textit{$k$-assignment} $L$ of a hypergraph $\hyh$ is a mapping from $V(\hyh)$ to the power set of $\N$ such that $|L(v)|=k$ for each vertex $v$ in $\hyh$.
For any $k$-assignment $L$, 
an \textit{$L$-coloring} of $\hyh$ is 
a proper coloring $\theta$ of $\hyh$ 
such that $\theta(v)\in L(v)$ for each $v\in V(\hyh)$.
Denote by $P(\hyh, L)$ the number of $L$-colorings of $\hyh$.
The \textit{list-color function} $P_l(\hyh,k)$ of a hypergraph $\hyh$,
introduced by Kostochka and Sidorenko~\cite{Kosto},
is defined to be the minimum value of $P(\hyh, L)$'s among all $k$-assignments $L$ of $\hyh$, for each $k\in\N$.
For the development of list-color functions, the readers can refer to~\cite{Dong22, Kosto, Thomassen,wang17}. 

\subsection{Results}

In 1992, Kostochka and Sidorenko~\cite{Kosto} raised the question of how similar 
the list-color function of a graph is to its chromatic polynomial.
Very soon, Donner~\cite{Donner} managed to show that for any graph $G$, $P(G,k)=P_l(G,k)$ holds when $k$ is sufficiently large, and in 2009, Thomassen~\cite{Thomassen} proved that $k> |V(G)|^{10}$ is a sufficient condition for $P(G,k)$ and $P_l(G,k)$ being equal. Later, Wang, Qian and Yan~\cite{wang17} improved this condition to $k\ge\frac{m-1}
	{\log(1+\sqrt{2})}\approx 1.1346(m-1)$,
	where $m=|E(G)|$.
The authors of this article~\cite{Dong22} recently obtained the conclusion
that  for any integer $k\ge m-1\ge 3$
and any $k$-assignment $L$ of $G$,
$P(G,L)-P(G,k)$ is bounded below by 
\equ{eq1-1}
{
\left ( (k-m+1)k^{n-3}+
	\frac{(m-1)(m-3)(k-m+3)k^{n-5}}{24}
	\right )
	\sum_{uv\in E(G)}|L(u)\setminus L(v)|,
}
where $n=|V(G)|$, 
implying that  $P_l(G,k)=P(G,k)$ holds.

In this paper, we focus on 
comparing 
the list-color function of a uniform hypergraph 
with its chromatic polynomial.
The best known result on this topic 
is the one below due to Wang, Qian and Yan in~\cite{wang20}.

\begin{theo}[\cite{wang20}]\label{wanghy}
For any $r$-uniform hypergraph $\hyh$ with $m$ edges,
where $r\ge 3$,
$P_l(\hyh,k)=P(\hyh,k)$ holds whenever 
$k\ge 1.1346(m-1)$.
\end{theo}

Now let $\hyh$ be an $r$-uniform hypergraph with $m$ edges, 
where $r\ge 2$. 
In Section~\ref{sec2},
we will establish a lower bound
for $P(\hyh, L)-P(\hyh,k)$
for any $k$-assignment $L$ of $\hyh$. 
Applying this result, 
we will show 
 in Section~\ref{sec3-1} 
that  when $m\le 4$, 
$P_l(\hyh,k)=P(\hyh,k)$ holds
for all $k\ge m-1$ (see 
Corollary~\ref{ssize}), and 
in Section~\ref{sec4}
that when $m\ge 5$, the following conclusion holds.

\begin{theo}\label{th4-1}
Let $\hyh$ be any $r$-uniform and  connected hypergraph with $n$ vertices and $m$  edges,
where $r\ge 2$ and $m\ge 5$.
For any $k$-assignment $L$ of $\hyh$,
if $k\ge m-1$, then 
\aln{eq4-1}
{
P(\hyh,L)-P(\hyh,k)
	\ge 
	(k-m+1)k^{n-r-1}
	\sum_{e\in E(\hyh)}
	\left ( 
	k-\left|\bigcap\limits_{v\in e}L(v)\right|
	\right ).
}
Therefore $P_l(\hyh,k)=P(\hyh,k)$ holds for any integer $k\ge m-1$.
\end{theo}

Theorem~\ref{th4-1}
follows directly from 
Propositions~\ref{pp1-5} and \ref{pr-new}, 
which will be established in 
Sections~\ref{sec2}
and \ref{sec4} respectively.

\section{A lower bound of $P(\hyh,L)-P(\hyh,k)$ 
	\label{sec2}}

In this section, we assume that $\hyh=(V,E)$ is a hypergraph with $n$ vertices and $m$ edges, $\eta$ is a bijection from $E$ to $[m]$, and $L$ is a $k$-assignment of $\hyh$. 
We shall establish a lower bound 
for $P(\hyh,L)-P(\hyh,k)$ in  Proposition~\ref{pp1-5}.

\begin{figure}[h]
	\centering 
	\includegraphics[width=12cm]	{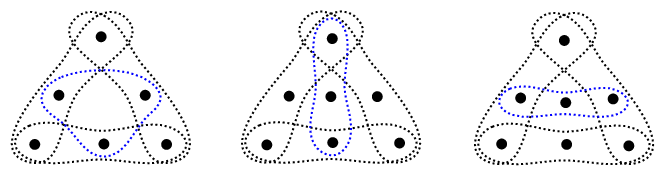}
	
	{}\hfill \hspace{0.5cm} (a) $\hyh_1$ \hspace{2.8 cm} 
	(b) $\hyh_2$ \hspace{2.8 cm}  (c) $\hyh_3$  \hfill {}    	
	\caption{Three $3$-uniform hypergraphs of size four}
	\label{fig1}
\end{figure}

Similar to Whitney's expression  for 
the chromatic polynomial of a graph in~\cite{Whitney1932},
the following expression for 
$P(\hyh,k)$ can also be obtained
 by applying the principle of inclusion-exclusion:
\aln{BC}
{
	P(\hyh,k)=\sum_{A\subseteq E}(-1)^{|A|}k^{c(A)}.
}
Recall that the idea of Whitney's broken-cycle theorem for the chromatic polynomial of a graph 
is to reduce as many pairs of edge sets as possible in the summation 
of $(-1)^{|A|}k^{c(A)}$.
Following the same idea, a broken-cycle theorem for proper colorings of hypergraphs has been established by Trinks in~\cite{Trinks14}. 
A \textit{$\delta$-cycle} in $\hyh$ is defined to be a minimal edge set $F$ in $\hyh$  such that $e\subseteq V(F\setminus \{e\})$ for every $e\in F$.
Clearly, in Figure~\ref{fig1}, 
the edge set of $\hyh_1$ is a $\delta$-cycle, while those of 
$\hyh_2$ and $\hyh_3$ are not.
If $\hyh$ is a graph, then a $\delta$-cycle in $\hyh$ is actually 
a cycle.
It is noteworthy that $\delta$-cycles may never exist in possibly many hypergraphs.

For each $\delta$-cycle $C$ in $\hyh$,
if $e$ is the edge in $C$ with $\eta(e)\le \eta(e')$ for all $e'\in C$,
then $C\setminus \{e\}$ is called a 
\textit{broken-$\delta$-cycle} of $\hyh$
with respect to $\eta$.
Let $\nbroken(\hyh)$ be the set of  subsets of $E$ containing no broken-$\delta$-cycles. 
Trinks~\cite{Trinks14} showed that  
the summation at the right-hand side of (\ref{BC}) can be restricted to the sets in $\nbroken(\hyh)$, 
as stated below:
\equ{thhyBC}
{
	P(\hyh,k)=\sum_{A\in \nbroken(\hyh)}(-1)^{|A|}k^{c(A)}.
}
For any $A\in \nbroken(\hyh)$, 
let $\beta(A,L)
=\prod\limits_{i=1}^{c(A)}\beta(\hyh_i,L)
$,
where 
$\hyh_1,\dots,\hyh_{c(A)}$ are the components of  $\hyh\langle A\rangle$
and $\beta(\hyh_i,L)=\big|\bigcap
\limits_{v\in V(\hyh_i)}L(v)\big|$.
Wang, Qian and Yan~\cite{wang20} deduced a similar result for $P(\hyh,L)$:
\equ{thhyBCL}
{
	P(\hyh,L)
	=\sum_{A\in \nbroken(\hyh)}(-1)^{|A|} \beta(A,L).
}
Then, 
by (\ref{thhyBC}) and (\ref{thhyBCL}), we have
\aln{eq3-1}
{
	P(\hyh,L)-P(\hyh,k)=\sum_{A\in \nbroken(\hyh)}(-1)^{|A|}\left(\beta(A,L)-k^{c(A)}\right).
}

A lower  bound of $\beta(A,L)-k^{c(A)}$ 
was provided by Wang, Qian and Yan~\cite{wang20}.
For any $e\in E$,
let $\alpha(e,L)=k-\left|\bigcap\limits_{v\in e}L(v)\right|$.

\begin{pro}[\cite{wang20}]
\label{pro3-1}
For any $k$-assignment $L$ of $\hyh$
	and any $A\subseteq E$,
	$$
	\beta(A,L)-k^{c(A)}\ge -k^{c(A)-1}\sum_{e\in A}\alpha(e,L).
	$$
\end{pro} 

An upper bound for $\beta(A,L)-k^{c(A)}$ will be provided 
in Proposition~\ref{pro3-2}
for any non-empty $A\subseteq E$ by applying the following lemma in \cite{Dong22}.

\begin{lem}
	[\cite{Dong22}]\label{lem14}
	Let $d_1,d_2, \dots, d_s$ be any non-negative real numbers, and $q_1,q_2,\dots,q_s$ be any positive real numbers, where $s\ge 1$. If $x\ge \max_{1\le i\le s}d_i$, then
	$$
	(x-d_1)(x-d_2)\cdots(x-d_s)\le x^s-\frac{x^{s-1}}{q_1+q_2+\cdots+q_s}\sum_{i=1}^sq_id_i.
	$$
\end{lem}

\prop{pro3-2}
{
	Let $L$ be a $k$-assignment of $\hyh$. Then, for any non-empty subset $A$ of $E$,
	$$
	\beta(A,L)-k^{c(A)}\le -\frac{k^{c(A)-1}}{|A|}\sum_{e\in A}\alpha(e,L).
	$$
}

\proof
Let $\hyh_1,\dots,\hyh_{c(A)}$ be the components of 
$\hyh\langle A\rangle$, where $E(\hyh_i)\ne \emptyset$ for all $i=1,2,\dots,q$ and $E(\hyh_i)=\emptyset$ for all $i=q+1,\dots,c(A)$. 
For any $i\in[q]$, 
by the definition of 
$\beta(\hyh_i,L)$ and $\alpha(e,L)$,
\aln{in3-2}
{
	\beta(\hyh_i,L)
	&=
	\big | \bigcap_{v\in V(\hyh_i)} L(v)
	\big |
	\le \min_{e\in E(\hyh_i)} 
	\big | \bigcap_{v\in e} L(v)
	\big |
	\nonumber \\
	&=k-\max_{e\in E(\hyh_i)}\alpha(e,L)\le k-\frac{1}{|E(\hyh_i)|}\sum_{e\in E(\hyh_i)}\alpha(e,L).
}
Then, by (\ref{in3-2}) and Lemma~\ref{lem14}, 
\equ{eq3-2}
{
	\prod_{i=1}^q\beta(\hyh_i,L)
	\le \prod_{i=1}^{q}\left(k-\frac{1}{|E(\hyh_i)|}\sum_{e\in E(\hyh_i)}\alpha(e,L)\right)
	\le  k^q-\frac{k^{q-1}}{|A|}\sum_{e\in A}\alpha(e,L).
} 
Consequently, the result follows as $\beta(\hyh_i,L)=k$ for all $i=q+1,\dots,c(A)$.
\proofend

For any $i\in[m]$ and $e\in E$, let
$\nbroken_i(\hyh)$ be the set of $A\in \nbroken(\hyh)$ with $|A|=i$,
$\nbroken(\hyh,e)$ be the set of $A\in \nbroken(\hyh)$ with $e\in A$,
and 
 $\nbroken_i(\hyh,e)
 =\nbroken_i(\hyh)\cap \nbroken(\hyh,e)$. 
 Obviously, $|\nbroken_i(\hyh,e)|\le \binom{m-1}{i-1}$. 

For  any $A\in\nbroken_i(\hyh)$,  
	if $i=1$, then $c(A)=n-r+1$;
	otherwise, 
	$e\not\subseteq V(A\setminus \{e\})$
	holds for some $e\in A$. 
	Applying this fact, 
	an upper bound for $c(A)$ was established in~\cite{wang20}.
	
	\begin{lem}[\cite{wang20}]
		\label{wanghyc}
		Let $\hyh$ be an $r$-uniform hypergraph with $n$ vertices, where $r\ge 2$.
		For any $A\in \nbroken_i(\hyh)$,
		where $i\ge 1$,  $c(A)\le n-r-i+2$.
	\end{lem}

For any edge $e$ in $\hyh$, 
let $E_{r-1}(e)$ be the set of edges $e'\in E$ with $|e'\cap e|=r-1$,
and 
let 
$\nbroken^*_2(\hyh, e)
=\{ \{e,e'\}\in \nbroken_2(\hyh, e):   e'\in E_{r-1}(e)\}$.
Clearly, $A\in \nbroken^*_2(\hyh, e)$
if and only if $A\in \nbroken_2(\hyh, e)$
with $c(A)=n-r$.
Now we present the main result of this section.

\prop{pp1-5}
{
	Let $\hyh$ be an $r$-uniform hypergraph with $n$ vertices and $m$ edges, where $r\ge 2$,  and $\eta$ be a bijection from $E$ to $[m]$.
 Then, for 	any $k$-assignment $L$ of $\hyh$, 
	\eqn{eq3-10}
	{
	P(\hyh, L)-P(\hyh,k)
	\ge 
k^{n-r} 
\sum_{e\in E}
\left (\alpha(e, L)
F_{\eta}(\hyh, e,k)\right ),
	}
where 
\eqn{eq3-11}
{
	F_{\eta}(\hyh,e,k)
	:&=&1-\frac{|\nbroken^*_2(\hyh, e)|}{k}
	-\frac{|\nbroken_2(\hyh,e)\setminus \nbroken^*_2(\hyh, e)|}{k^2}
	\nonumber \\
& &+\sum_{2\le i\le (m+1)/2 \atop 
		A\in \nbroken_{2i-1}(\hyh,e)}
	\frac {k^{c(A)-1-(n-r)}}{2i-1}
	-  	\sum_{2\le i\le m/2 \atop 
		A\in \nbroken_{2i}(\hyh,e)}
	k^{c(A)-1-(n-r)}.
}
}

\proof
By  (\ref{eq3-1}) and Propositions~\ref{pro3-1} and~\ref{pro3-2},
we have 
\aln{eq3-4}
{
	P(\hyh,L)-P(\hyh,k)
	&=
	-\sum_{A\in \nbroken(\hyh)\atop |A|~\text{odd}}\left(\beta(A,L)-k^{c(A)}\right)
	+\sum_{A\in \nbroken(\hyh)\atop |A|~\text{even}}\left(\beta(A,L)-k^{c(A)}\right)
	\nonumber\\
	&\ge
	\sum_{A\in \nbroken(\hyh)\atop |A|~\text{odd}}
	\left(\frac{k^{c(A)-1}}{|A|}\sum_{e\in A}\alpha(e,L)\right)+
	\sum_{A\in \nbroken(\hyh)\atop |A|~\text{even}}
	\left(-k^{c(A)-1}\sum_{e\in A}\alpha(e,L)\right)
	\nonumber\\
	&=
	\sum_{e\in E}
	\left (\alpha(e,L)
	\left(\sum_{A\in \nbroken(\hyh,e)\atop |A|~\text{odd}}\frac{k^{c(A)-1}}{|A|}
	-\sum_{A\in \nbroken(\hyh,e)\atop |A|~\text{even}}k^{c(A)-1}
	\right)
	\right ).
}
Note that $c(\{e\})=n-r+1$.
For $A\in \nbroken_2(\hyh,e)$, 
we have  $c(A)\le n-r$, and 
$c(A)=n-r$ if and only if 
$A\in \nbroken^*_2(\hyh,e)$.
	Thus, the result follows from 
	the 
	definition of $F_{\eta}(\hyh,e,k)$.
	\proofend

Clearly, $|\nbroken^*_2(\hyh,e)|\le |E_{r-1}(e)|$ and 
$|\nbroken_2(\hyh,e)\setminus \nbroken^*_2(\hyh,e)|\le m-1-|E_{r-1}(e)|$.
By the expression of (\ref{eq3-11}),
$F_{\eta}(\hyh,e,k)$ is independent  of $L$.

\section{$m\le 4$
	\label{sec3-1}
}

In this section, we consider
any $r$-uniform hypergraph
$\hyh=(V,E)$ with $m=|E|\le 4$.

\lemm{le-n0}
{ 
	Let $\hyh=(V,E)$ be a connected $r$-uniform hypergraph with $m$ edges, where $r\ge 2$,
	and 
	$\eta$ be any 
	bijection 
	from $E$ to $[m]$.
	Then, 	for each edge $e$ in $\hyh$, 
	$F_{\eta}(\hyh,e,k)\ge 0$ holds
	whenever $1\le m\le 4$ and $k\ge m-1$. 
}

\proof Let  $e\in E$. 
If $m=1$, then $F_{\eta}(\hyh, e,k)=1$
for any $k\ge 0$.
Now assume that $2\le m\le 4$. 
We first show that 
\equ{le-n01}
{
	F_{\eta}(\hyh, e,k)
	\ge 1-\frac{|E_{r-1}(e)|}{k}
	-\frac{m-1-|E_{r-1}(e)|}{k^2}.
}
If $m=2$, or $3\le m\le 4$ 
and $\hyh$ has $\delta$-cycles, then
$\nbroken_4(\hyh)=\emptyset$ 
and  (\ref{le-n01}) follows from (\ref{eq3-11}).

If $m=4$ and $\hyh$ has no $\delta$-cycle,
then 
$\nbroken_4(\hyh)=\{E\}$
and 
$\nbroken_3(\hyh)=
\{S\subseteq E: |S|=3\}$. 
Note that $c(E)=1$ as $\hyh$ is connected.
Then by (\ref{eq3-11}), 
\eqn{n2-02}
{
	F_{\eta}(\hyh,e,k)
	&\ge & 1-\frac{|E_{r-1}(e)|}{k} -\frac{m-1-|E_{r-1}(e)|}{k^2} 
	+\frac 13 \sum_{A\in \nbroken_3(\hyh,e)}k^{c(A)-1-(n-r)}
	-k^{c(E)-1-(n-r)}
	\nonumber \\
	&\ge&  1-\frac{|E_{r-1}(e)|}{k} -\frac{m-1-|E_{r-1}(e)|}{k^2} 
	+\frac 13 \times {3\choose 2}k^{1-1-(n-r)}
	-k^{1-1-(n-r)}
	\nonumber \\
	&=& 1-\frac{|E_{r-1}(e)|}{k} -\frac{m-1-|E_{r-1}(e)|}{k^2} .
}
Hence (\ref{le-n01}) holds. 
Obviously,   (\ref{le-n01}) implies that 
$F_{\eta}(\hyh,e,k)\ge 1-\frac{m-1}{k}
\ge 0$ whenever $k\ge m-1$. 
The result holds.
\proofend 

Due to Proposition~\ref{pp1-5}, Lemma~\ref{le-n0} and the fact that both the chromatic polynomial and the list-color function factorize over the components of a hypergraph, 
we have the following conclusion.

\corr{ssize}
{
	Let $\hyh=(V,E)$ be an $r$-uniform hypergraph, where $r\ge 2$ and $m=|E|\le 4$.
	Then $P_l(\hyh,k)=P(\hyh,k)$ holds 
	for $k\ge m-1$.
}

\section{$m\ge 5$
\label{sec4} }

In this section, we 
assume that $\hyh=(V,E)$ is an $r$-uniform hypergraph
with $r\ge 2$ and $m=|E|\ge 5$, 
and $\eta$ is a bijection from $E$ to $[m]$.
Throughout this section, 
we assume that 
$e_0$ is a fixed edge in $\hyh$.

\subsection{A lower bound of $F_{\eta}(\hyh, e_0,k)$ for any edge $e_0$
	\label{nsec3}
}

Recall that 
$E_{r-1}(e_0)$ is the set of edges $e\in E$ with $|e\cap e_0|=r-1$.
For any $i\ge 0$, let $\setd_i$ be the set of $D\subseteq E_{r-1}(e_0)$ such that $|D|=i$  and $\{e_0\}\cup D\in \nbroken(\hyh)$.
By  definition,  $\setd_{i}=\emptyset$ 
for all $i>|E_{r-1}(e_0)|$.
Our purpose in this subsection is
to show in Proposition~\ref{npro-3} 
that 
$F_{\eta}(\hyh, e_0,k)\ge 
f_{m-1,k, |E_{r-1}(e_0)|}(|\setd_0|, |\setd_1|, \dots, |\setd_{m}|)$,
where 
$f_{M,k,N}(z_0, z_1, \dots , z_{M+1})$
is a function depending on non-negative numbers 
$M,k,N, z_0, z_1,  \dots , z_{M+1}$  only,
where $M,k$ and $N$ are integers.

In the next subsection, we shall show that 
$f_{M,k,N}(z_0, z_1, \dots , z_{M+1})
\ge \frac{k-M}{k}$ holds
under the conditions that 
$k\ge M\ge N\ge z_0+z_1+\cdots+z_{M+1}$,
$0\le z_i\le {N\choose i}$ 
for all $0\le i\le M+1$ 
and 
$kz_{2i}\ge (2i+1)z_{2i+1}$ 
for all  
$0\le i\le M/2$.
Then, it follows that 
$F_{\eta}(\hyh, e_0,k)\ge \frac{k-M}{k}$.

We first establish some conclusions on 
$\setd_i$'s.

\prop{npro-1}
{
$\setd_i=0$ for all $i\ge m$, and for any $0\le i\le m-1$, 
	\begin{enumerate} 
\item 
$|\nbroken_{i+1}(\hyh,e_0)|\ge |\setd_i|$; 
\item $ |\setd_i|\le {|E_{r-1}(e_0)\choose i}$;  and 
\item 
$(|E_{r-1}(e_0)|-i)|\setd_i|\ge (i+1)|\setd_{i+1}|$.
\end{enumerate} 
}

\proof  By definition, $\setd_i=\emptyset$ 
when $i\ge m$. 
(i) and (ii) follow directly from definition.
(iii) follows from the facts that 
\vspace{-0.3 cm}
\begin{enumerate}
	\item [(a)] for any $D\in \setd_{i+1}$ and $e\in D$, 
we have $D\setminus \{e\}\in \setd_i$,
and 
\item [(b)]
for any $D'\in \setd_{i}$,
there are at most $|E_{r-1}(e_0)|-i$
edges $e'\in E\setminus D'$ 
such that $D'\cup \{e'\}\in \setd_{i+1}$. 
\end{enumerate}
\proofend

\prop{npro-2}
{
	For any $i$ with $0\le i\le m-1$, 
	\equ{b2}
	{
		|\nbroken_{i+1}(\hyh,e_0)|
		\le \sum_{0\le t\le i}
		{m-1-|E_{r-1}(e_0)|\choose i-t}
		|\setd_t|.
	}
}

\proof For any $A\in \nbroken_{i+1}(\hyh, e_0)$, 
$A$ can be partitioned into three subsets
of $\nbroken(\hyh)$:
$\{e_0\}$, $A\cap E_{r-1}(e_0)$
and $A\setminus (\{e_0\}\cup E_{r-1}(e_0))$. 

For any integer $t$ with $0\le t\le i$, 
by definition, 
the number of subsets $A_1$ of $E_{r-1}(e_0)$ 
such that $A_1\cup \{e_0\}\in \nbroken_{t+1}(\hyh,e_0)$ 
is equal to $|\setd_t|$, 
and the number of subsets $A_2$ 
of $E\setminus (\{e_0\}\cup E_{r-1}(e_0))$ with $A_2\in \nbroken_{i-t}(\hyh)$ is at most 
${|E|-1-|E_{r-1}(e_0)|\choose i-t}$.
Thus, (\ref{b2}) holds. 
\proofend

Let $U$ be the set of vertices 
$u\in V\setminus e_0$ such that 
$u\in e\subseteq  e_0\cup \{u\}$
for some edge $e$ in $\hyh$.
Clearly, $U$ is the union of $e\setminus e_0$ over all edges $e\in E_{r-1}(e_0)$.
Assume that $U=\{u_1,\dots,u_s\}$, where $s=|U|$.
For any $i\in [s]$, 
let $E(e_0,u_i)$ be the set of edges $e\in E$ 
with $u_i\in e\subseteq e_0\cup \{u_i\}$.
Thus, $E_{r-1}(e_0)$ is partitioned 
into subsets 
$E(e_0,u_1), E(e_0,u_2), \dots, E(e_0,u_s)$.

\lemm{nle-3}
{For any $D\subseteq E_{r-1}(e_0)$,
	if 
	$\{e_0\}\cup D\in \nbroken(\hyh)$,
	then
	$|D\cap E(e_0,u_i)|\le 1$ 
	for each $i\in [s]$,
	and thus  $\setd_i=\emptyset$ 
	for $i>s$.
}

\proof 
Assume that 
$\{e_0\}\cup D\in \nbroken(\hyh)$, where $D\subseteq E_{r-1}(e)$. 
Observe that for any $e,e'\in E(e_0,u_i)$, 
where $i\in [s]$, 
$\{e_0,e,e'\}$ forms a $\delta$-cycle.
Thus,
$|D\cap E(e_0,u_i)|\le 1$ for each $i\in [s]$,
implying that $\setd_i=\emptyset$ 
for all $i>s$.
\proofend

We  are now going to establish a lower bound for $F_{\eta}(\hyh,e_0,k)$ in terms of the following function: 
\eqn{npro-3-1}
{
f_{M,k,N}(Z)
=	1-\frac{N}{k}
	-\frac{M-N}{k^2}
	+\sum_{1\le i\le M/2}
	\frac {k^{-2i}}{2i+1}
	z_{2i}
	-  	\sum_{1\le i\le M/2\atop
		0\le t\le 2i+1}
	k^{-2i-1}
	{M-N\choose 2i+1-t}z_{t}.
}
where $Z$ is the sequence $(z_0,z_1,z_2,\dots, z_{M+1})$.

\prop{npro-3}
{
$F_{\eta}(\hyh,e_0,k)\ge f_{m-1,k,|E_{r-1}(e_0)|}(Z)$,
where 
$Z=(|\setd_0|, |\setd_1|, |\setd_2|,
\dots, |\setd_{m}|)$.
}

\proof We will apply Proposition~\ref{pp1-5}
to prove this result. 
Recall that $|\nbroken^*_2(\hyh,e)|\le |E_{r-1}(e)|$ and 
$|\nbroken_2(\hyh,e)\setminus \nbroken^*_2(\hyh,e)|\le m-1-|E_{r-1}(e)|$. Thus
\equ{nle-29}
{
1-\frac{|\nbroken^*_2(\hyh, e)|}{k}
-\frac{|\nbroken_2(\hyh,e)\setminus \nbroken^*_2(\hyh, e)|}{k^2}
\ge 1-\frac{|E_{r-1}(e_0)|}k - 
\frac{m-1-|E_{r-1}(e_0)|}{k^2}.
}

Let $D$ be any member in $ \setd_{2i}$, where  $i\ge 1$.
By Lemma~\ref{nle-3}, $i\le s/2$ and 
$|D\cap E(e_0,u_j)|\le 1$ 
for each $j\in [s]$.
It follows that 
$c(D\cup\{e_0\})=n-r-|D|+1$. 
By definition, $D\cup \{e_0\}\in  \nbroken_{2i+1}(\hyh,e_0)$. 
Thus, 
\eqn{nle-30}
{\sum_{2\le i\le (m+1)/2 \atop A\in \nbroken_{2i-1}(\hyh,e_0)}\frac {k^{c(A)-1-(n-r)}}{2i-1} 
&=&\sum_{1\le i\le (m-1)/2 \atop A\in \nbroken_{2i+1}(\hyh,e_0)}\frac {k^{c(A)-1-(n-r)}}{2i+1}
	\nonumber \\
&\ge&
\sum_{1\le i\le (m-1)/2 \atop 
D\in \setd_{2i}
}
\frac {k^{c(D\cup \{e_0\})-1-(n-r)}}{2i+1}
	\nonumber \\
&=&\sum_{1\le i\le (m-1)/2}\frac {k^{-2i}}{2i+1}	|\setd_{2i}|.
}

By Lemma~\ref{wanghyc}, 
$c(A)\le  n-r-(2i+2)+2$ for each 
$A\in \nbroken_{2i+2}(\hyh, e_0)$. 
Proposition~\ref{npro-2} then implies that 
\eqn{nle-31}
{\sum_{2\le i\le m/2 \atop A\in \nbroken_{2i}(\hyh,e_0)}k^{c(A)-1-(n-r)}
	&=&\sum_{1\le i\le (m-2)/2 \atop A\in \nbroken_{2i+2}(\hyh,e_0)}
	k^{c(A)-1-(n-r)}
	\nonumber\\
	&\le&	\sum_{1\le i\le (m-2)/2\atop
		0\le t\le 2i+1 }
	k^{-2i-1}
	{m-1-|E_{r-1}(e_0)|\choose 2i+1-t}
	|\setd_t|.
}
The result then follows from Proposition~\ref{pp1-5} and (\ref{nle-29}), (\ref{nle-30}) and  (\ref{nle-31}).
\proofend

\subsection{$f_{M,k,N}(Z)
	\ge \frac{k-M}{k}$ under certain conditions
\label{sec5}
}

In this section, 
we shall prove that 
$f_{M,k,N}(Z)
\ge \frac{k-M}{k}$ 
holds for positive integers $k,M,N$
with $k\ge M\ge N$ and 
any sequence  $Z=(z_0,z_1,z_2,\dots, z_{M+1})$ satisfying 
the conditions that 
$0\le z_i\le {N\choose i}$ 
for all $0\le i\le M+1$ 
and 
$kz_{2i}\ge (2i+1)z_{2i+1}$ 
for all  
$0\le i\le M/2$.
We first show that $f_{M,k,N}(z_0,z_1,\dots, z_M)$ 
is a decreasing function on $N$
for any fixed numbers $k,M$ and 
$z_0,z_1,\dots,z_M$ 
under certain conditions.

\lemm{nle-10}
{
	Assume that $M,N$ and $k$ are 
	integers with 
	$M\ge 4$ and $k\ge M>N\ge 0$.
	If $0\le z_i\le {N\choose i}$ for each 
	$i$ with $0\le i\le M+1$, 
	then 
	$f_{M,k,N}(z_0,z_1,\dots, z_{M+1})\ge f_{M,k,N+1}(z_0,z_1,\dots, z_{M+1})$.
}

\proof Let $Z=(z_0,z_1,\dots, z_{M+1})$. 
By the definition of 
$f_{M,k,N}(Z)$ in (\ref{npro-3-1}), 
\eqn{nle-10-e1}
{
	f_{M,k,N}(Z)- f_{M,k,N+1}(Z)
	&=&\frac 1{k}-\frac{1}{k^2}
	-\sum_{1\le i\le M/2\atop
		0\le t\le 2i+1 }
	k^{-2i-1}
	{M-1-N\choose 2i-t} z_t
	\nonumber \\ &\ge &
	\frac 1{k}-\frac{1}{k^2}
	-\sum_{1\le i\le M/2\atop
		0\le t\le 2i+1 }
	k^{-2i-1}
	{M-1-N\choose 2i-t}
	{N\choose t} 
	\nonumber \\ &=&
	\frac 1{k}-\frac{1}{k^2}
	-\frac 1k \sum_{1\le i\le M/2}
	k^{-2i}
	{M-1\choose 2i}
	\nonumber \\ &\ge &
	\frac 1{k}-\frac{1}{k^2}
	-\frac 1{k} 
	\left (-1+
	\sum_{i\ge 0}
	k^{-2i}\frac{(M-1)^{2i}}{(2i)!}\right )
	\nonumber \\ &= &
	\frac 2{k}-\frac{1}{k^2}
	-\frac 1{2k} 
	\left ( 
	\exp((M-1)/k)+\exp((1-M)/k)
	\right) .
}
Define 
\equ{nle-10-e2}
{
	\phi(x,M)
	=2x - x^2 - \frac x2
	(\exp((M-1)x) + \exp((1-M)x)).
}
Since $k\ge M$, 
by (\ref{nle-10-e1}), 
it remains to show that 
$\phi(x,M)\ge 0$ for 
$0\le x\le \frac1{M}\le \frac1{M-1}$. 

It can be verified  that $\exp((M-1)x) + \exp((1-M)x)$ is an increasing function of $x$ in $[0, \infty)$ for any fixed $M\ge 4$.
Thus, for $0\le x\le \frac 1{M-1}\le \frac 13$, we have 
\equ{nle-10-e3}
{
	\phi(x,M)
	\ge 2x - x^2 - \frac x2
	(\exp(1) + \exp(-1))
	\ge 0.4569x-x^2\ge 0.
}
Hence Lemma~\ref{nle-10} follows from 
(\ref{nle-10-e1}) and (\ref{nle-10-e3}).
\proofend

Now we are going to prove the following conclusion. 

\prop{lm-new}
{
		Assume that $k, M,  N$ are integers with  
$M\ge 4$ and $k\ge M\ge N\ge 0$.
	If 
	$0\le z_i\le {N\choose i}$ 
	for all $0\le i\le M+1$ 
	and 
	$kz_{2i}\ge (2i+1)z_{2i+1}$ 
	for all  
	$0\le i\le M/2$,
	then 
	$f_{N,k,M}(z_0,z_1,\dots, z_M)\ge
	\frac{k-M}{k}$.
}

\proof 
By Lemma~\ref{nle-10} and the given 
conditions, 
\eqn{spec-e0}
{
f_{N,k,M}(Z)\ge f_{M,k,M}(Z)
&=& 
1-\frac{M}{k}
+\sum_{1\le i\le M/2} \frac {k^{-2i}}{2i+1}
z_{2i}
-  	\sum_{1\le i\le (M-1)/2}
k^{-2i-1}z_{2i+1}
\nonumber \\
&\ge &
1-\frac{M}{k}
+\sum_{1\le i\le M/2}  k^{-2i-1}\left ( 
\frac {k\cdot z_{2i}}{2i+1}
-  z_{2i+1}\right )
\nonumber \\
&\ge & 
1-\frac{M}{k}.
}
\proofend

\begin{pro}
	\label{pr-new}
	For any  $r$-uniform  hypergraph $\hyh=(V,E)$, where $r\ge 2$ 
	and  $m=|E|\ge 5$,
	if  $k\ge m-1$,  then 
	$F_{\eta}(\hyh,e_0,k)\ge \frac{k-(m-1)}{k}$ holds for each $e_0\in E$.
\end{pro} 

\proof  
The result follows from Propositions~\ref{pp1-5},~\ref{npro-3} and~\ref{lm-new}
by taking $M=m-1$ and $N=|E_{r-1}(e_0)|$. 
The conditions 
that 
$z_i\le {N\choose i}$ for all $0\le i\le M$
and  $kz_i\ge (i+1)z_{i+1}$ 
for all  $0\le i\le M-1$
follow from
Proposition~\ref{npro-1} (ii) and (iii)
respectively. 
\proofend 

Propositions~\ref{pp1-5} and \ref{pr-new} imply 
Theorem~\ref{th4-1} directly.

\section{Further study}

Recall that for any hypergraph 
$\hyh=(V,E)$, the \textit{chromatic number} (or resp., \textit{list chromatic number}) of $\hyh$, denoted by $\chi(\hyh)$ (or resp., $\chi_l(\hyh)$), is the minimum integer $q$ such that $P(\hyh,q)>0$ (or resp., $P_l(\hyh,q)>0$) holds.
Obviously, 
$|V|\ge \chi_l(\hyh)\ge \chi(\hyh)$. 
We wonder if 
$P_l(\hyh,k)=P(\hyh,k)$ holds 
for all hypergraphs $\hyh=(V,E)$
and integers $k\ge |V|$.

\section*{Acknowledgement}

This research is supported by the Ministry of Education,
Singapore, under its Academic Research Tier 1 (RG19/22). Any opinions,
findings and conclusions or recommendations expressed in this
material are those of the author(s) and do not reflect the views of the
Ministry of Education, Singapore).

\end{document}